\documentclass{amsart}
\usepackage{amssymb}

\theoremstyle{definition}

\theoremstyle{remark}

\numberwithin{equation}{section}

\subjclass[2000]{Primary 20D08}

\def\squareforqed{\hbox{\rlap{$\sqcap$}$\sqcup$}}
\def\qed{\ifmmode\squareforqed\else{\unskip\nobreak\hfil
\penalty50\hskip1em\null\nobreak\hfil\squareforqed
\parfillskip=0pt\finalhyphendemerits=0\endgraf}\fi\medskip}

\newcommand{\udot}{{}^{\textstyle .}}

\newcommand{\PSL}{\mathrm{PSL}}

\newcommand{\Aut}{\mathrm{Aut}\,}
\newcommand{\Sz}{\mathrm{Sz}}

\title{Maximal subgroups of sporadic groups}
\date{8th April, 2016}
\author{Robert A. Wilson}

\address{School of Mathematical Sciences\\
Queen Mary University of London\\
London E1 4NS\\U.K.}
\email{R.A.Wilson@qmul.ac.uk}
\begin{document}
\maketitle

\begin{abstract}
A systematic study of maximal subgroups of the sporadic simple groups began in the 1960s.
The work is now almost complete, only a few cases in the Monster remaining outstanding.
We give a survey of results obtained, and methods used,
over the past $50$ years, for the classification of maximal
subgroups of sporadic simple groups, and their automorphism groups.
\end{abstract}

\section{Introduction}
The subtitle of the `Atlas of Finite Groups' \cite{Atlas}
is `Maximal Subgroups and Ordinary Characters
for Simple Groups'. These two aspects of the study of finite simple groups remain at the
forefront of research today.
The Atlas was dedicated to collecting facts, not to providing proofs. It contains an
extensive bibliography, but not citations at the point of use, making it difficult for the 
casual reader to track down proofs. In the ensuing 30 years, moreover, the landscape has
changed dramatically, both with the appearance of new proofs in the literature, and
with the ability of modern computer algebra systems to recompute much of the data
in the twinkling of an eye.

As far as maximal subgroups are concerned,
shortly before the publication of the Atlas it became clear
that the maximal subgroup project should be extended to almost simple groups.
The reason for this is that it is not possible to deduce the maximal subgroups of an
almost simple group directly from the maximal subgroups of the corresponding
simple group. This was made clear by the examples described in \cite{MaxAut}, especially
perhaps the maximal subgroup $S_5$ of $M_{12}{:}2$, which is 
neither the normalizer of a maximal subgroup of $M_{12}$, nor the normalizer
of the intersection of two non-conjugate maximal subgroups of $M_{12}$.
 
The results on maximal subgroups for all the classical groups
in the Atlas, as well as exceptional groups of types ${}^2B_2$, ${}^2G_2$, $G_2$
and ${}^3D_4$, are proved, and in many instances corrected, 
in the recent book by Bray, Holt and Roney-Dougal
\cite{BHRD}. This leaves the sporadic groups and five exceptional groups,
${}^2F_4(2)'$, $F_4(2)$, $E_6(2)$, ${}^2E_6(2)$, $E_7(2)$ and $E_8(2)$.
Of the latter, completeness of the list of maximal subgroups was claimed only for
${}^2F_4(2)'$, and a reference is given to \cite{RuTits}, although a correction
is noted in \cite {ABC}.
Subsequently, three cases have been completed and published, namely
$F_4(2)$ by Norton and Wilson \cite{F42}, $E_6(2)$ by Kleidman and Wilson
\cite{E62}, and $E_7(2)$ by Ballantyne, Bates and Rowley \cite{E72}. The case
$E_8(2)$ is still not complete, while the proof for ${}^2E_6(2)$
has not been published. I have re-calculated the latter case, and can confirm
that the list in the Atlas is complete.

Turning now to the sporadic groups and their automorphism groups, just $7$ of the
$26$ cases were unfinished at the time of publication of the Atlas, namely the three
Fischer groups, $J_4$, the Thompson group, the Monster and the Baby Monster.
All except the Monster have now been completed, and lists of maximal
subgroups for the simple groups (although not for the automorphism groups)
are given in \cite{TFSG} (but note two known errors: the subgroup $O_8^+(3){:}S_4$
is missing from the list of maximal subgroups of the Baby Monster on page 261,
and the subgroup $5^2{:}2A_5$ of $Co_1$ is given wrongly as $5^2{:}4A_5$ on page 211). 
In this paper we survey the results and main methods,
and try to provide full references where the details may be found.

\section{Methods}
\subsection{Basic strategy}
Suppose $S$ is a sporadic simple group, and $S\le G\le\Aut(S)$, and suppose
$M\ne S$ is a maximal subgroup of $G$. Then $M\cap S\ne 1$, and if $K$ is a
minimal characteristic subgroup of $M\cap S$, then $M=N_G(K)$.
Since minimal characteristic subgroups of finite groups are direct products of
isomorphic simple groups, this leads to a $3$-step process:
\begin{enumerate}
\item determine the characteristically simple subgroups $K$ of $S$, up to conjugacy;
\item determine the normalizer $N_G(K)$ in each case;
\item decide maximality in $G$ in each case.
\end{enumerate}
Usually Step 1 is by far the most difficult.

There is a fundamental difference between the case when $K$ is an elementary abelian
$p$-group (the so-called \emph{$p$-local} case), and the case when $M\cap S$ has
no non-trivial abelian normal subgroup (the so-called \emph{non-local} case, in which 
$K$ is necessarily non-abelian).
Occasionally $M\cap S$ may have more than one minimal characteristic subgroup, and
it may be convenient to dispose of these cases first to avoid duplication of effort.
Indeed, even if $M\cap S$ has a unique minimal characteristic subgroup, $M$ 
itself may not,
and such cases may also be treated separately.

\subsection {Local subgroups}
The techniques of local analysis developed for the classification project for finite
simple groups are powerful enough to go quite a long way towards the determination
of the maximal $p$-local subgroups. One of the first things to be worked out
when a new sporadic group was found was the list of conjugacy classes of elements of
prime order, and the corresponding centralizers. This provides a convenient basis
for an inductive classification of all $p$-local subgroups.

One starts by finding the conjugacy classes of elements of order $p$ in
$N_S(\langle x\rangle)/\langle x \rangle$, where $x$ has order $p$. This gives a list of
subgroups of order $p^2$, and one next determines which of these groups are
conjugate, at the same time computing their centralizers, and normalizers.
Discarding the cyclic groups of order $p^2$, one can then proceed to the next
level, classifying the elements of order $p$ in
$N_S(\langle x,y\rangle)/\langle x,y \rangle$, where
$\langle x,y\rangle$ is elementary abelian of order $p^2$. And so on.
If $G$ has a suitable representation, this process can even be automated, as was
done for example by Greg Butler \cite{ButlerHe} for the $2$-local
subgroups of the Held group.

Once the Sylow $p$-subgroup gets large, however, this brute-force approach
becomes cumbersome, and various refinements are required, especially for $p=2$.
For example, if $S$ has a double cover, then the squaring map on an elementary
abelian $2$-group in $S$ lifts to a quadratic form, which we may assume is either
zero or non-singular. Indeed, quadratic or symplectic forms can sometimes be constructed in
other cases where there are two classes of involutions, even if there is no double cover.
For example, in my PhD thesis I showed that in the Fischer group
$Fi_{22}$, an elementary abelian $2$-group which is $2B$-pure
supports an invariant symplectic form, defined to be $1$ when the $3$-transposition factors
of the two $2B$-elements fail to commute with each other. In \cite{Filoc}
every elementary abelian $2$-subgroup of $Fi_{24}'$ is shown to support an invariant
symplectic form.

Similar ideas were used for example by Kleidman and Wilson \cite{J4max}
in $J_4$ and by Meierfrankenfeld and Shpectorov \cite{MSh,Meier} in
the Monster. More generally, if there is
more than one class of elements of order $p$, there is scope for a creative
division into the cases that need to be considered.

\subsection{Existence of non-local subgroups}
In the search for, and classification of, non-local subgroups, 
the methods will vary depending on whether one is looking for a new simple group or not.
I shall assume the latter, and therefore use CFSG when necessary. First consider the
problem of proving \emph{existence} of non-local subgroups. This is surprisingly hard.
For example, the embedding of $M_{12}$ in $M_{24}$ was unknown for many
years, until proved by Frobenius. When $J_1$ was discovered \cite{J1}, it was conjectured
that it contained a subgroup isomorphic to $L_2(11)$, but to confirm this
required substantial (for those days) computational assistance. Later, Livingstone \cite{Liv}
constructed $J_1$ as the automorphism group of a graph on $266$ vertices,
thereby providing an alternative proof. More recently, the existence of subgroups
of the Monster isomorphic to $L_2(71)$
 \cite{A5subs}, 
$L_2(59)$ 
\cite{L259}
and $L_2(41)$ \cite{L241} has been
verified only once, again using substantial computational resources.

There are some theoretical methods that can be used, but they have limited applicability.
These include the so-called Brauer trick, in which a subgroup is constructed from an
amalgam of two subgroups $H$ and $K$ with specified intersection. If there is a
representation in which the dimensions of the fixed-point spaces $U,V,W$
of $H$, $K$ and
$H\cap K$ satisfy $$\dim U + \dim V > \dim W,$$ then the group
generated by $H$ and $K$ has a fixed point, so is a proper subgroup.
More generally, the amalgam method can sometimes be used to identify the
subgroup, even when there is no fixed point.

\subsection{Non-existence of non-local subgroups}
Usually, the hard part of the classification of non-local maximal subgroups is the proof that
there are no more than the ones that have been constructed. Non-local subgroups
with non-simple socle are relatively easy to classify, as every composition factor of the socle
has to lie in the centralizer of an element of prime order $p$, for some $p\ge 5$, which
generally reduces the possibilities to a manageable list. So the problem reduces
essentially to classifying the non-abelian simple subgroups up to conjugacy.

A reasonable first step is to try to classify them up to isomorphism, using
Lagrange's Theorem together with character restriction and
knowledge of the $p$-local subgroups to
eliminate as many cases as possible. In larger cases, these methods are rarely
sufficient, and it will be necessary to carry along a number of `possible'  isomorphism
types of simple subgroups for more detailed investigation. For example, Griess \cite{FG}
was unable, even with the referee's help, to decide if $J_1$ is a subgroup of the
Monster. This was eventually decided in the negative \cite{J1M}. The last case in
the sporadic groups, namely the question whether $Sz(8)$ is contained
in the Monster, was decided (in the negative) only very recently \cite{Sz8M}.

More advanced techniques which we shall discuss in more detail include use of
structure constants, restriction of Brauer characters (using \cite{ABC}), 
computation of
cohomology, detailed local analysis,
and extending normalizers.
The first and last of these techniques often reduce to
an exhaustive computational search. There is a general tendency
to work inductively from small simple groups to large, as exemplified most 
clearly in Norton's work on the Monster \cite{Anatomy1}, using the
normalizer extension method.

\section{Historical survey}
\subsection{Livingstone and his students: structure constants}
It would appear that the first systematic attempts to classify the maximal subgroups
of sporadic simple groups were undertaken by Donald Livingstone and his students
in the decade following the discovery of $J_1$. While the case of $J_1$ itself was
straightforward \cite{J1}, 
apart from the difficulty of establishing the existence of a subgroup
$L_2(11)$, the same was not necessarily true for the five previously known
sporadic groups, namely the Mathieu groups. 

Chang Choi \cite{Choi1,Choi2}
began with $M_{24}$, and completed his thesis in 1968.
A particularly troublesome case was the classification of transitive imprimitive subgroups,
largely because there was a previously unsuspected maximal subgroup of this type,
isomorphic to $L_2(7)$.
This subgroup was apparently first
found by Robert Curtis \cite{CurtisM24}, who went on to
provide a new proof for the list of maximal subgroups of $M_{24}$,
using his newly-discovered
Miracle Octad Generator \cite{MOG}. This work forms part of his
thesis, completed in 1972 under John Conway. 
Another proof was published by Rudy List \cite{list}.
As for the four smaller Mathieu groups,
it is not clear that proofs of completeness of their lists of maximal subgroups have
even been published. Certainly the literature search conducted while preparing the 
Atlas did not throw up any such references. Nevertheless, these days such proofs
would be regarded as graduate student exercises, and do not present great difficulty.

Two more students of Livingstone, both of whom completed their theses in 1970,
were Spyros Magliveras, who determined the maximal subgroups of
the Higman--Sims group \cite{Magliveras}, and Larry Finkelstein
who did the same for 
the McLaughlin group and Conway's third group $Co_3$ \cite{Finkelstein}.
Finkelstein then collaborated with Arunas Rudvalis to deal with the Janko groups
$J_2$ \cite{FRJ2} and $J_3$ \cite{FRJ3}. 

Among the various techniques they used was the method of
structure constants. Given three conjugacy classes $C_1,C_2,C_3$ in $G$ there
is a well known character formula which counts the number of ways a fixed element
$z\in C_3$ can be written as the product $z=xy$ 
of elements $x\in C_1$ and $y\in C_2$.
If the orders of $x,y,z$ are $p,q,r$ respectively, then in certain cases the
isomorphism type of $\langle x,y,z\rangle$ is determined:
$$
\begin{array}{|ccc|c|}
p&q&r&\langle x,y,z\rangle\cr
\hline
2&2&n& D_{2n}\cr
2&3&3& A_4\cr
2&3&4&S_4\cr
2&3&5& A_5\cr\hline
\end{array}
$$
This is particularly useful for classifying subgroups isomorphic to $A_5$.

Of note here, however, is the use that was put to calculating structure constants for
other triples of integers, especially $(p,q,r)=(2,3,7)$. Since, for example, $L_2(7)$
is generated by elements $x$ of order $2$ and $y$ of order $3$ with $xy$ of order $7$,
calculation of the structure constants of type $(2,3,7)$ gives an upper bound
on the number of subgroups isomorphic to $L_2(7)$, and in some cases this bound
is actually met.

\subsection{Students of Conway: lattice methods}
By the mid-1970s the age of discovery of sporadic simple groups was over, and there
were some $14$ new groups whose maximal subgroups were waiting to be determined.
Up to this point, every case considered had a relatively small permutation
representation: the largest was $Co_3$, on $276$ points. The next generation of groups
needed thousands of points, and demanded new techniques of investigation.

The work of Robert Curtis on maximal subgroups of $M_{24}$ was only part of
his thesis, completed in 1972, under John Conway.
He also considered in depth the classification of subgroups
of the largest Conway group \cite{Slattices,Co1loc}. At that time he did not
envisage a complete determination of the maximal subgroups, which had to wait
another decade. Of particular interest here is his approach to studying subgroups
which fix interesting sublattices. This is an extension of Conway's discovery that
the Higman--Sims and McLaughlin groups, originally constructed as permutation
groups, are also essentially sublattice stabilizers in the Conway group. Curtis 
discovered several interesting subgroups of $Co_1$ by these methods.
He also completely classified the maximal $p$-local subgroups \cite{Co1loc}.
(The case $p=3$ had been done earlier by Mikdashi, a student of Livingstone,
but not published except in his PhD thesis, dated 1971. See also \cite{Co13loc} for a 
correction to \cite{Co1loc} for $p=3$.)

Two more of Conway's students worked on subgroups of sporadic groups
in the 1970s, and worked out a good deal of the subgroup structure,
although apparently without aiming at a full classification
of the maximal subgroups. 
Simon Norton wrote his thesis in 1975 on the group now known as
the Harada--Norton group, among other topics. This thesis does not
claim to determine the maximal subgroups completely, 
although he did this later. This work was however not published until
a decade later, when I collaborated with him
to re-work the determination of the maximal subgroups of the Harada--Norton group
\cite{HN} and its automorphism group.
Gerard Enright wrote his thesis in 1977 on the subgroup structure of the Fischer groups
$Fi_{22}$ and $Fi_{23}$, in which his main result was a classification of the
subgroups generated by transpositions \cite{Enright}.

Norton went on to do a huge amount of work on the subgroup structure of the Monster,
and by extension, the Baby Monster and the Harada--Norton group, as well as
$Fi_{24}'$ and other groups involved in the Monster. His big idea was to create a
table of `Monstralizer pairs', that is, pairs of subgroups $H$ and $K$ such that
$H=C(K)$ and $K=C(H)$. The published version of this table \cite{Anatomy1}
contains all cases where $H$ or $K$ has order divisible by a prime $p$ with $p\ge 11$,
but in unpublished work he went further than this.

In 1979 I became a student of Conway, and
worked first on the Suzuki group (1782 points) \cite{WilSuz}, and then
the Rudvalis group (4060 points) \cite{RuTits}. The methods were not much different
from the methods of Livingstone's students, 
except perhaps in a greater emphasis on using properties of
lattices on which the groups (or, more accurately, their covering groups) act. For example,
the 6-fold cover of the Suzuki group acts on the complex Leech lattice, of dimension $12$.
If one can show by character theory that a subgroup isomorphic to $K$ has a fixed 
vector in the lattice, then the same is true in the lattice reduced modulo $3$. But there are
only two orbits of non-zero vectors in the latter, and the vector stabilizers are already known.
Thus the difficult cases to classify are the irreducible subgroups, in which case
exhaustive computer searches were carried out. This included construction of the subgroups
$A_7$, $L_2(25)$ and $L_3(3)$. Similarly, in the Rudvalis group, the subgroups
$L_3(3)$, $L_2(13)$ and $A_8$
were found by computer searches. The subgroups $U_3(5)$ and $L_2(29)$ had
been found much earlier  by Kiang-Chuen Young (a student of
John McKay in Canada, PhD thesis 1974), using similar computational methods.

The lattice method really showed its power, however, in the case of $Co_2$, where I used
a $23$-dimensional sublattice of the Leech lattice
\cite{Co2}. The fact that $23$ is prime was
particularly helpful, as it implied that most proper subgroups were reducible, and very
few of those did not have a fixed vector. Indeed, character restriction alone is sufficient to 
show that every proper non-abelian simple subgroup of $Co_2$ fixes a vector in
the $23$-dimensional representation. Automorphisms of order $2$ swap two fixed vectors and therefore fix their sum, so the only slight difficulty arises for groups which have automorphisms
of order $3$, in this case $L_2(8)$, $L_3(4)$ and $U_3(5)$. 

It seemed natural then to try to apply the same
methods to $Co_1$, acting on the Leech lattice itself
\cite{Co1}. In this case the so-called Suzuki chain of subgroups provides a long list
of (quasi-)simple subgroups which act fixed point freely on the Leech lattice,
and the degree $24$ allows a number of other cases as well, such as $L_2(11)$,
$L_2(23)$ and $L_2(25)$. The alternating groups, for example, were 
classified by first using the
structure constants of type $(2,4,5)$ to limit the number of possibilities
for $A_6$, and then inductively constructing
$A_n$ from $A_{n-2}$ and $A_{n-3}\times 3$ intersecting in $A_{n-3}$.
Only the final three cases, $L_2(11)$,
$L_2(23)$ and $L_2(25)$ required computer calculations
to complete.
At this stage it became clear that the lattice method could also usefully
be applied to
the cases considered by Magliveras and Finkelstein, namely the Higman--Sims group,
the McLaughlin group and Conway's third group. This provided independent proofs
for the complete
classification of non-local subgroups in these cases \cite{vecstab}.  

Although my thesis, in common with earlier work, 
only dealt with the simple groups, and not their automorphism
groups, it became clear fairly soon afterwards that it was necessary to deal with the
latter case also. In \cite{MaxAut} I went through all the relevant cases up to that
point, and adjusted the proofs to include the automorphism groups as well, that is,
the groups $M_{12}{:}2$, $M_{22}{:}2$, $J_2{:}2$, $J_3{:}2$, 
$HS{:}2$, $McL{:}2$, $Suz{:}2$, $He{:}2$ and $HN{:}2$. 

\subsection{Global input}
By the early 1980s, then, 
the number of target groups had been reduced to single figures, and a number of
people around the world were attacking them. 
Greg Butler in Australia (a student of John Cannon) dealt with the Held group,
using serious computational methods, for example in the systematic enumeration
of $2$-local subgroups.
Satoshi Yoshiara in Japan had written his
thesis (unfortunately in Japanese) on the Suzuki group, and went on to determine
the maximal subgroups of the O'Nan group \cite{YoshiON}. Independently, in Moscow,
Ivanov, Tsaranov and Shpectorov \cite{ITSON} did the same, as did I \cite{WilON}.
Our methods were quite different, and it was reassuring to find that we all obtained the
same answer. The most difficult part was constructing subgroups isomorphic to $A_7$,
$L_2(31)$, and $M_{11}$, for which I used computation, whereas the other authors
above used detailed geometrical methods.

It was also becoming evident that the problems were getting harder, and from this point on
it was rare for the maximal subgroup problem for a single group to be completely
solved in a single paper. The Lyons group was considered by Andrew Woldar \cite{WolLy}
(a student of Ron Solomon) in his PhD thesis in 1984, in which he also
conjectured the existence of a $111$-dimensional representation over the field of order $5$.
This representation was also conjectured by Meyer and Neutsch, and
constructed by Richard Parker \cite{MNP}, and was used extensively in my work
to complete the determination of the maximal subgroups \cite{WilLy1,WilLy2}.

\subsection{Students of Wilson: hard-core computation}
It had long been clear that computational methods were a necessary part of
the maximal subgroup project for sporadic groups. Theoretical methods were
simply not powerful enough to probe the structure of these unique objects
in sufficient detail.
In terms of checking the results, it may be felt that these computer searches represent
a weak point. Nowadays, however, computer algebra systems such as MAGMA and GAP
are sufficiently well developed that it should be possible with relatively little effort,
and insignificant amounts of computer time, to reproduce these results robustly.
As far as I am aware, this effort has not yet been made, but it is surely time now
to do so.

The story of the original hard-core computation 
really starts with Peter Kleidman, who was a student of Martin Liebeck.
When Liebeck moved from Cambridge to London, Peter worked much of the time with me.
We first finished the determination of the maximal subgroups of $Fi_{22}$,
and its automorphism group $Fi_{22}{:}2$,
which I had been unable to complete on my own \cite{WilF22, KWF22}. 
Then we attacked $J_4$ \cite{J4max}, largely
because it had a small representation, of dimension $112$ over the field of order $2$.
We found two new maximal subgroups $U_3(3)$ 
and $M_{22}{:}2$. Independently,
Wolfgang Lempken \cite{Lempken} classified the maximal $p$-local subgroups of $J_4$,
but his methods were insufficient to complete the determination of the non-local
subgroups. Kleidman and I were then joined by Richard Parker in a project to
classify the maximal subgroups of $Fi_{23}$.

My first official PhD student was Steve Linton, whose thesis, dated 1990, was on maximal
subgroups of the Thompson sporadic simple group $Th$ \cite{LintonTh}
and the Fischer group
$Fi_{24}'$ \cite{F24} and its automorphism group. 
The $p$-local subgroups of $Fi_{24}'$ had been classified in \cite{Filoc}.
Linton's work on $Th$ also built on my earlier work \cite{someTh}, which reduced
the problem to classifying subgroups isomorphic to
$L_2(19)$, $A_6$, $L_2(7)$, $L_3(3)$ and $U_3(3)$. 
In particular, he discovered new subgroups $L_2(19){:}2$ and $L_3(3)$. Much of the work
was done using computations in the $248$-dimensional representation.
A major theme of the computations was a systematic enumeration of pairs of
elements $(x,y)$ with $x$, $y$ and $xy$ in specified conjugacy classes.

My next student to 
work on maximal subgroups of sporadic groups was Petra (Beth) Holmes, whose PhD thesis
on `Computing in the Monster' dates from 2002. But the Monster really deserves
a section to itself, as perhaps does the Baby Monster.

\subsection{The Baby Monster}
Simon Norton had already done a lot of work on the maximal subgroups of the
Baby Monster, but most of this has not been published.
In \cite{someB} I had classified the $p$-local subgroups for $p$ odd, and provided
some basic information about non-local subgroups. The $2$-local subgroups were
classified by Meierfrankenfeld and Shpectorov \cite{MSh,Meier}.
In \cite{moreB} about half of the cases of the non-local subgroup problem were
dealt with, by theoretical methods. 
This relied heavily on Norton's work on subgroups of the Monster.

Then in \cite{BMcon} I constructed generators for the Baby Monster as
$4370\times 4370$ matrices over the field of
order $2$, which paved the way for a complete determination of the maximal
subgroups. Nevertheless, the computations were not straightforward, as at that
time a single matrix multiplication took around $20$ minutes.
The first result \cite{newB} was the construction of previously unknown
maximal subgroups $L_2(31)$ and $L_2(49)\udot 2$.
The project was completed in \cite{maxB}, although
the promised follow-up paper containing details never appeared.

One technique employed here, which became even more important in the Monster,
was to classify subgroups generated by two copies of $A_5$ intersecting in $D_{10}$.
Since the subgroups $A_5$ and $D_{10}$
can both be counted using the structure constants, it is relatively straighforward to
enumerate all the cases, and test for isomorphism with any desired group.
Other cases were more like the case of $L_2(17)$, which can be generated by subgroups
$17{:}8$ and $D_{16}$ intersecting in a cyclic group of order $8$. In this case,
I found representatives for the three classes of $17{:}8$ in the Baby Monster,
and then found the normalizers of the three cyclic groups of order $8$, using
standard methods \cite{Bray}
for finding involution centralizers, repeatedly.  
\subsection{The Monster}
The Monster is a special case because of its enormous size, and a large number 
of papers (at least $15$) 
have been written on various aspects of its subgroup structure. To date $44$ conjugacy
classes of maximal subgroups are known, and the possibilities for presently unknown
maximal subgroups are severely restricted.
The $p$-local maximal 
subgroups for odd $p$ were determined in \cite{oddlocals}, and
for $p=2$ in \cite{MSh,Meier}.

The first serious attack on the non-local subgroups was carried out by Simon Norton,
whose results are reported in \cite{Anatomy1}, although without detailed proofs.
From a classification of subgroups isomorphic to $A_5$, obtained in part from the
calculation of structure constants, he obtained a complete
classification of simple subgroups containing an $A_5$ with $5A$-elements. 
The method is illustrated by the case $A_6$, which can be generated by two copies
of $A_5$, intersecting in $D_{10}$ (or $S_4$: both methods are useful). Hence the
centralizer of any $A_6$ is the intersection of two copies of the $A_5$-centralizer
inside the $D_{10}$-centralizer. The latter is either the Harada--Norton group or
an involution centralizer therein, so the calculation reduces to computing double cosets
of certain subgroups of $HN$.

This work was extended in \cite{Anatomy2}, which classified simple groups into three
categories: those which were definitely in the Monster, those which were definitely not,
and those  for which we could not decide at that time. These lists are not quite correct,
as $L_2(41)$ was put in the `definitely not' category, whereas we now know that
$L_2(41)$ is in fact a subgroup of the Monster \cite{L241}. The case $J_1$ was handled
in \cite{J1M}.

At that point, a computer construction of the Monster became available, based on the
$3$-local subgroup $3^{1+12}\udot2\udot Suz{:}2$. This was used by my student
Beth Holmes to investigate subgroup structure, but it soon became clear that we
needed an involution centralizer, so we built the Monster again  \cite{2loccon}, this time using
$2^{1+24}\udot Co_1$. The first result of this work was the discovery of a
new maximal subgroup $L_2(29){:}2$ \cite{L229}, followed closely by
$L_2(59)$ \cite{L259}. This was then extended to a systematic study of 
subgroups generated by two copies of $A_5$ with $5B$-elements \cite{A5subs}, which
turned up new maximal subgroups $L_2(71)$ and $L_2(19){:}2$.
Later, Holmes \cite{S4subs} classified subgroups isomorphic to $S_4$, and used this to classify
subgroups isomorphic to $U_3(3)$, $L_3(3)$, $L_2(17)$, and $L_2(7)$.
More recent computations were used to classify subgroups $L_2(41)$
\cite{L241}, $L_2(27)$ \cite{L227}, and $L_2(13)$ containing $13B$-elements 
\cite{L213B}. A largely theoretical, but delicate, proof that the Monster does not
contain $Sz(8)$ is given in \cite{Sz8M}.

At the time of writing, 
the published results on non-local subgroups of the Monster include complete
classifications of maximal subgroups with simple socle, in all
cases except when the socle is $L_2(8)$, $L_2(13)$, $L_2(16)$, $U_3(4)$ or $U_3(8)$.
In more recent work, not yet published, I have eliminated the case
$U_3(8)$ theoretically, and the cases $U_3(4)$ and $L_2(8)$ computationally.
The remaining two cases do not appear to present significantly greater difficulty,
and should be completed before long. 
\section{Results}
There are very few known errors in the lists of maximal subgroups in the Atlas:
\begin{itemize}
\item in $J_3$ and $J_3{:}2$ the shape of the Sylow $3$-subgroup is given as
$3^2.(3\times 3^2)$, which should be $3^2.3^{1+2}$;
\item $L_2(17){:}2$ is wrongly included as a subgroup of $Fi_{23}$;
\item in $Co_1$, the two groups described as $N(3C^2)$, of shapes
$3^2.[2.3^6].2A_4$ and $3^2.[2^3.3^4].2A_4$, either do not exist or are not
maximal. Also, the subgroup described as $5^2{:}4A_5$ is actually
$5^2{:}2A_5$. 
\item in the Monster, the subgroup described as $(A_7\times (A_5\times A_5).4).2$
should be described as $(A_7\times (A_5\times A_5).2^2).2$.
\end{itemize}
Additional information obtained since the publication of the Atlas includes the
following:
\begin{itemize}
\item the lists for $Fi_{22}$ and $Fi_{22}{:}2$ are complete;
\item for $Th$, all listed subgroups exist and are maximal, and the list becomes complete after
adding $L_3(3)$;
\item for $Fi_{23}$, the list becomes complete after deleting $L_2(17){:}2$ and
adding $L_2(23)$;
\item for $J_4$, the list becomes complete after adding $M_{22}{:}2$ and $U_3(3)$;
\item for $Fi_{24}$ the list is complete, and for $Fi_{24}'$ becomes complete after
adding two classes of $L_2(13){:}2$.
\item for the Baby Monster, all $p$-local subgroups which are listed without an overgroup
are maximal, and the list of maximal subgroups becomes complete on adding the
following eight classes of non-local subgroups: $(S_6\times L_3(4){:}2){:}2$,
$(S_6\times S_6).4$, $L_2(49)\udot 2$, $L_2(31)$, $M_{11}$, $L_3(3)$, $L_2(17){:}2$,
and $L_2(11){:}2$.
\item for the Monster, there is one $7$-local maximal subgroup to be added, of
shape $7^2{:}SL_2(7)$, and five non-local maximal subgroups:
$L_2(71)$, $L_2(59)$, $L_2(41)$, $L_2(29){:}2$, and $L_2(19){:}2$.
Any further maximal subgroup has socle $L_2(13)$ or $L_2(16)$.
\end{itemize}

We conclude with tables of results which give an update on the Atlas. In each case, the
maximal subgroups are listed in decreasing order of order. To save space, the cases of
two conjugacy classes of maximal subgroups of $S$ which are fused in $\Aut(S)$ are
denoted by the annotation `(two)'. The listing of two maximal subgroups of $M_{12}{:}2$
of shape $L_2(11){:}2$ is not a mistake: there is no automorphism fusing these classes
of subgroups, and they are fundamentally different.

$$
\begin{array}{c}
M_{11}\cr\hline A_6\udot 2\cr L_2(11)\cr 3^2{:}SD_{16}\cr S_5\cr 2\udot S_4\cr\hline
\end{array}
\;
\begin{array}{c}
J_1\cr\hline L_2(11)\cr 2^3{:}7{:}3\cr 2\times A_5\cr 19{:}6\cr 11{:}10\cr
D_6\times D_{10}\cr 7{:}6\cr\hline
\end{array}
\;
\begin{array}{c}
M_{22}\cr\hline L_3(4)\cr 2^4{:}A_6\cr A_7 \mbox{ (two)}\cr 2^4{:}S_5\cr
2^3{:}L_3(2)\cr A_6\udot 2\cr L_2(11)\cr\hline
\end{array}
\;
\begin{array}{c}
M_{22}{:}2\cr\hline M_{22}\cr L_3(4){:}2_2\cr 2^4{:}S_6\cr 2^5{:}S_5\cr
2^3{:}L_3(2)\times 2\cr A_6\udot 2^2\cr L_2(11){:}2\cr\hline
\end{array}
\;
\begin{array}{c}
M_{23}\cr\hline M_{22}\cr L_3(4){:}2_2\cr 2^4{:}A_7\cr A_8\cr M_{11}\cr
2^4{:}(3\times A_5){:}2\cr 23{:}11\cr\hline
\end{array}
\;
\begin{array}{c}
M_{12}\cr\hline
M_{11} \mbox{ (two)}\cr A_6\udot 2^2 \mbox{ (two)}
\cr L_2(11)\cr 3^2{:}2S_4 \mbox{ (two)}\cr
 2\times S_5\cr 2^{1+4}S_3\cr 4^2{:}D_{12}\cr A_4\times S_3\cr\hline
\end{array}
$$
$$
\begin{array}{c}
M_{12}{:}2\cr\hline
M_{12}\cr L_2(11){:}2\cr L_2(11){:}2\cr  (2^2\times A_5){:}2\cr
2^{1+4}.D_{12}\cr 4^2.D_{12}.2\cr 3^{1+2}{:}D_8\cr S_4\times S_3\cr S_5\cr\hline
\end{array}
\;
\begin{array}{c}
J_3\cr\hline L_2(16){:}2\cr L_2(19) \mbox{ (two)}\cr 2^4{:}(3\times A_5)\cr L_2(17)\cr
(3\times A_6){:}2_2\cr 3^2.3^{1+2}{:}8\cr 2^{1+4}{:}A_5\cr 2^{2+4}{:}(3\times S_3)\cr\hline
\end{array}
\;
\begin{array}{c}
J_3{:}2\cr\hline J_3\cr L_2(16){:}4\cr 2^4{:}(3\times A_5).2\cr L_2(17)\times 2\cr
(3\times M_{10}){:}2\cr 3^2.3^{1+2}{:}8.2\cr 2^{1+4}S_5\cr 2^{2+4}{:}(S_3\times S_3)\cr 19{:}18\cr \hline
\end{array}
\;
\begin{array}{c}
M_{24}\cr\hline M_{23}\cr M_{22}{:}2\cr 2^4{:}A_8\cr M_{12}{:}2\cr 2^6{:}3\udot S_6\cr
L_3(4){:}S_3\cr 2^6{:}(L_3(2)\times S_3)\cr L_2(23)\cr L_2(7)\cr\hline
\end{array}
\;
\begin{array}{c}
J_2\cr\hline U_3(3)\cr 3\udot PGL_2(9)\cr 2^{1+4}{:}A_5\cr 2^{1+4}{:}(3\times S_3)\cr
A_4\times A_5\cr A_5\times D_{10}\cr L_3(2){:}2\cr 5^2{:}D_{12}\cr A_5\cr\hline
\end{array}
$$

$$
\begin{array}{c}
Ly\cr\hline G_2(5)\cr 3\udot McL{:}2\cr 5^3\udot L_3(5)\cr 2\udot A_{11}\cr
5^{1+4}{:}4S_6\cr 3^5{:}(2\times M_{11})\cr 3^{2+4}{:}2A_5.D_8\cr 67{:}22\cr
37{:}18\cr\hline
\end{array}
\;
\begin{array}{c}
O\mbox{'}N\cr\hline L_3(7){:}2  \mbox{ (two)}\cr J_1\cr 4_2\udot L_3(4){:}2_1\cr
(3^2{:}4\times A_6)\udot 2\cr 3^4{:}2^{1+4}D_{10}\cr L_2(31) \mbox{ (two)}\cr
4^3\udot L_3(2)\cr M_{11} \mbox{ (two)}\cr A_7 \mbox{ (two)}\cr\hline
\end{array}
\;
\begin{array}{c}
O\mbox{'}N{:}2\cr\hline O\mbox{'}N\cr J_1\times 2\cr
4_2\udot L_3(4)\udot 2^2\cr (3^2{:}4\times A_6)\udot 2^2\cr 3^4{:}2^{1+4}D_{10}.2\cr
4^3\udot(L_3(2)\times 2)\cr  7^{1+2}{:}(3\times D_{16})\cr
31{:}30\cr  L_2(7){:}2\cr PGL_2(9)\cr\hline
\end{array}
\;
\begin{array}{c}
J_2{:}2\cr\hline J_2\cr G_2(2)\cr 3\udot A_6\udot 2^2\cr 2^{1+4}S_5\cr
2^{2+4}.(S_3\times S_3)\cr (A_4\times A_5){:}2\cr (A_5\times D_{10})\udot 2\cr
L_3(2){:}2\times 2\cr 5^2{:}(4\times S_3)\cr S_5\cr\hline
\end{array}
\;
\begin{array}{c}
HS\cr\hline M_{22}\cr U_3(5){:}2 \mbox{ (two)}\cr L_3(4){:}2_1\cr S_8\cr 2^4.S_6\cr
4^3{:}L_3(2)\cr M_{11} \mbox{ (two)}\cr 4\udot 2^4{:}S_5\cr 2\times A_6\udot 2^2\cr
5{:}4\times A_5\cr\hline
\end{array}
\;
$$

$$
\begin{array}{c}
HS{:}2\cr\hline HS\cr M_{22}{:}2\cr L_3(4){:}2^2\cr S_8\times 2\cr
2^5.S_6\cr 4^3{:}(L_3(2)\times 2)\cr 2^{1+6}S_5\cr (2\times A_6\udot 2^2).2\cr
 5^{1+2}{:}[2^5]\cr 5{:}4\times S_5\cr\hline
\end{array}
\;
\begin{array}{c}
McL\cr\hline U_4(3)\cr M_{22} \mbox{ (two)}\cr U_3(5)\cr 3^{1+4}{:}2S_5\cr 3^4{:}M_{10}\cr
L_3(4){:}2\cr 2\udot A_8\cr 2^4{:}A_7 \mbox{ (two)}\cr M_{11}\cr 5^{1+2}{:}3{:}8\cr\hline
\end{array}
\;
\begin{array}{c}
McL{:}2\cr\hline McL\cr U_4(3){:}2_3\cr U_3(5){:}2\cr 3^{1+4}{:}4S_5\cr
3^4{:}(M_{10}\times 2)\cr L_3(4){:}2^2\cr 2\udot S_8\cr 
M_{11}\times 2\cr 5^{1+2}{:}3{:}8.2\cr2^{2+4}{:}(S_3\times S_3)\cr\hline
\end{array}
\;
\begin{array}{c}
He\cr\hline S_4(4){:}2\cr 2^2\udot L_3(4).S_3\cr 2^6{:}3\udot S_6 \mbox{ (two)}\cr
2^{1+6}.L_3(2)\cr 7^2{:}2\udot L_2(7)\cr 3\udot S_7\cr 7^{1+2}{:}(S_3\times 3)\cr
S_4\times L_3(2)\cr 7{:}3\times L_3(2)\cr 5^2{:}4A_4\cr\hline
\end{array}
\;
\begin{array}{c}
Co_2\cr\hline U_6(2){:}2\cr 2^{10}{:}M_{22}{:}2\cr McL\cr 2^{1+8}{:}S_6(2)\cr HS{:}2\cr
(2^{1+6}\times 2^4)A_8\cr U_4(3).D_8\cr 2^{4+10}(S_5\times S_3)\cr M_{23}\cr
3^{1+4}{:}2^{1+4}S_5\cr 5^{1+2}{:}4S_4\cr\hline
\end{array}
$$

$$
\begin{array}{c}
He{:}2\cr\hline He\cr S_4(4){:}2\cr 2^2\udot L_3(4).D_{12}\cr
2^{1+6}.L_3(2).2\cr 7^2{:}2\udot L_2(7).2\cr
3\udot S_7\times 2\cr (S_5\times S_5){:}2\cr 2^{4+4}.(S_3\times S_3).2\cr 
7^{1+2}{:}(S_3\times 6)\cr S_4\times L_3(2){:}2\cr
7{:}6\times L_3(2)\cr 5^2{:}4S_4\cr\hline
\end{array}
\;
\begin{array}{c}
Fi_{22}\cr\hline 2\udot U_6(2)\cr O_7(3) \mbox{ (two)}\cr O_8^+(2){:}S_3\cr
2^{10}{:}M_{22}\cr 2^6{:}S_6(2)\cr (2\times 2^{1+8}{:}U_4(2)){:}2\cr
S_3\times U_4(3){:}2_2\cr {}^2F_4(2)'\cr 2^{5+8}{:}(S_3\times A_6)\cr
3^{1+6}{:}2^{3+4}{:}3^2{:}2\cr S_{10} \mbox{ (two)}\cr M_{12}\cr\hline
\end{array}
\;
\begin{array}{c}
Fi_{22}{:}2\cr\hline Fi_{22}\cr 2\udot U_6(2){:}2\cr 
O_8^+(2){:}S_3\times 2\cr 2^{10}{:}M_{22}{:}2\cr 2^7{:}S_6(2)\cr
(2\times 2^{1+8}{:}U_4(2){:}2){:}2\cr S_3\times U_4(3).2^2\cr {}^2F_4(2)\cr
2^{5+8}{:}(S_3\times S_6)\cr 3^5{:}(2\times U_4(2){:}2)\cr
3^{1+6}{:}2^{3+4}{:}3^2{:}2.2\cr G_2(3){:}2\cr M_{12}{:}2\cr\hline
\end{array}
\;
\begin{array}{c}
J_4\cr\hline 2^{11}{:}M_{24}\cr 2^{1+12}\udot 3\udot M_{22}{:}2\cr 2^{10}{:}L_5(2)\cr
2^{3+12}\udot(S_5\times L_3(2))\cr U_3(11){:}2\cr M_{22}{:}2\cr 
11^{1+2}{:}(5\times 2S_4)\cr L_2(32){:}5\cr L_2(23){:}2\cr U_3(3)\cr 29{:}28\cr 43{:}14\cr
37{:}12\cr\hline
\end{array}
\;
$$

$$
\begin{array}{c}
HN\cr\hline A_{12}\cr 2\udot HS.2\cr U_3(8){:}3\cr 2^{1+8}(A_5\times A_5).2\cr
(D_{10}\times U_3(5))\udot 2\cr 5^{1+4}{:}2^{1+4}.5.4\cr 2^6\udot U_4(2)\cr
(A_6\times A_6)\udot D_8\cr 2^{3+2+6}(3\times L_3(2))\cr 5^{2+1+2}{:}4A_5\cr
M_{12}{:}2 \mbox{ (two)}\cr 3^4{:}2(A_4\times A_4).4\cr 3^{1+4}{:}4A_5\cr\hline
\end{array}
\;
\begin{array}{c}
HN{:}2\cr\hline HN\cr S_{12}\cr 4\udot HS.2\cr U_3(8){:}6\cr 2^{1+8}(A_5\times A_5).2^2\cr
5{:}4 \times U_3(5){:}2\cr 5^{1+4}.2^{1+4}.5.4.2\cr 2^6\udot U_4(2).2\cr
(S_6\times S_6){:}2^2\cr 2^{3+2+6}(S_3\times L_3(2))\cr 5^{2+1+2}4S_5\cr
3^4{:}2(S_4\times S_4).2\cr 3^{1+4}{:}4S_5\cr\hline
\end{array}
\;
\begin{array}{c}
Co_3\cr\hline McL{:}2\cr HS\cr U_4(3){:}2^2\cr M_{23}\cr 3^5{:}(M_{11}\times 2)\cr
2\udot S_6(2)\cr U_3(5){:}S_3\cr 3^{1+4}{:}4S_6\cr 2^4\udot A_8\cr L_3(4){:}D_{12}\cr
2\times M_{12}\cr 2^2.[2^7.3^2].S_3\cr S_3\times L_2(8){:}3\cr A_4\times S_5\cr\hline
\end{array}
\;
\begin{array}{c}
Fi_{23}\cr\hline 2\udot Fi_{22}\cr O_8^+(3){:}S_3\cr 2^2\udot U_6(2).2\cr S_8(2)\cr
S_3\times O_7(3)\cr 2^{11}\udot M_{23}\cr 3^{1+8}.2^{1+6}.3^{1+2}.2S_4\cr
3^3.[3^7].(2\times L_3(3))\cr S_{12}\cr (2^2\times 2^{1+8}).(3\times U_4(2)).2\cr
2^{6+8}{:}(A_7\times S_3)\cr S_4\times S_6(2)\cr S_4(4){:}4\cr L_2(23)\cr\hline
\end{array}
\;
$$

$$
\begin{array}{c}
B \mbox{ (first part)}\cr\hline
2\udot {}^2E_6(2){:}2\cr 2^{1+22}\udot Co_2\cr Fi_{23}\cr 2^{9+16}.S_8(2)\cr
Th\cr (2^2\times F_4(2)){:}2\cr 2^{2+10+20}.(M_{22}{:}2\times S_3)\cr [2^{30}].L_5(2)\cr
S_3\times Fi_{22}{:}2\cr [2^{35}].(S_5\times L_3(2))\cr\hline
\end{array}
\;
\begin{array}{c}
B \mbox{ (second part)}\cr\hline
HN{:}2\cr O_8^+(3){:}S_4\cr 3^{1+8}.2^{1+6}.U_4(2).2\cr

(3^2{:}D_8\times U_4(3).2.2).2\cr 5{:}4 \times HS{:}2\cr S_4\times {}^2F_4(2)\cr
[3^{11}].(S_4\times 2S_4)\cr S_5\times M_{22}{:}2\cr (S_6\times L_3(4){:}2){:}2\cr
5^3\udot L_3(5)\cr\hline
\end{array}
\;
\begin{array}{c}
B \mbox{ (third part)}\cr\hline
5^{1+4}.2^{1+4}.A_5.4\cr (S_6\times S_6).4\cr 5^2{:}4S_4\times S_5\cr L_2(49)\udot 2_3\cr
L_2(31)\cr M_{11}\cr L_3(3)\cr L_2(17){:}2\cr L_2(11){:}2\cr 47{:}23\cr\hline
\end{array}
$$

$$
\begin{array}{c}
Ru\cr\hline {}^2F_4(2)\cr 2^6\udot G_2(2)\cr (2^2\times Sz(8)){:}3\cr 2^{3+8}{:}L_3(2)\cr
U_3(5){:}2\cr 2^{1+4+6}{:}S_5\cr L_2(25)\udot 2^2\cr A_8\cr L_2(29)\cr 5^2{:}4S_5\cr
3\udot A_6\udot 2^2\cr 5^{1+2}{:}[2^5]\cr L_2(13){:}2\cr A_6\udot 2^2\cr
5{:}4\times A_5\cr\hline
\end{array}
\;
\begin{array}{c}
Suz\cr\hline G_2(4)\cr 3_2\udot U_4(3){:}2\cr U_5(2)\cr 2^{1+6}\udot U_4(2)\cr
3^5{:}M_{11}\cr J_2{:}2\cr 2^{4+6}{:}3A_6\cr (A_4\times L_3(4)){:}2\cr
2^{2+8}{:}(A_5\times S_3)\cr M_{12}{:}2\cr 3^{2+4}{:}2(A_4\times 2^2).2\cr
(A_6\times A_5)\udot2\cr (3^2{:}4\times A_6)\udot 2\cr L_3(3){:}2 \mbox{ (two)}\cr
L_2(25)\cr A_7\cr\hline
\end{array}
\;
\begin{array}{c}
Suz{:}2\cr\hline Suz\cr G_2(4){:}2\cr 3_2\udot U_4(3){:}2^2\cr U_5(2){:}2\cr
2^{1+6}\udot U_4(2).2\cr 3^5{:}(M_{11}\times 2)\cr J_2{:}2\times 2\cr 2^{4+6}{:}3S_6\cr
(A_4\times L_3(4){:}2_3){:}2\cr 2^{2+8}{:}(S_5\times S_3)\cr M_{12}{:}2\times 2\cr
3^{2+4}{:}2(S_4\times D_8)\cr (PGL_2(9)\times A_5){:}2\cr (3^2{:}8\times A_6)\udot 2\cr
L_2(25){:}2\cr S_7\cr\hline
\end{array}
\;
\begin{array}{c}
Th\cr\hline {}^3D_4(2){:}3\cr 2^5\udot L_5(2)\cr 2^{1+8}\udot A_9\cr U_3(8){:}6\cr 
(3\times G_2(3)){:}2\cr 3.[3^8].2S_4\cr 3^2.[3^7].2S_4\cr 3^5{:}2S_6\cr 5^{1+2}{:}4S_4\cr
5^2{:}GL_2(5)\cr 7^2{:}(3\times 2S_4)\cr L_2(19){:}2\cr L_3(3)\cr M_{10}\cr 31{:}15\cr
S_5\cr\hline
\end{array}
\;$$

$$
\begin{array}{c}
Co_1\cr\hline Co_2\cr 3\udot Suz{:}2\cr 2^{11}{:}M_{24}\cr Co_3\cr
2^{1+8}\udot O_8^+(2)\cr U_6(2){:}S_3\cr (A_4\times G_2(4)){:}2\cr
2^{2+12}{:}(A_8\times S_3)\cr 2^{4+12}\udot (S_3\times 3S_6)\cr
3^2\udot U_4(3).D_8\cr 3^6{:}2M_{12}\cr (A_5\times J_2){:}2\cr
3^{1+4}{:}2\udot U_4(2){:}2\cr (A_6\times U_3(3)){:}2\cr 3^{3+4}{:}2(S_4\times S_4)\cr
A_9\times S_3\cr (A_7\times L_2(7)){:}2\cr (D_{10}\times (A_5\times A_5).2).2\cr
5^{1+2}{:}GL_2(5)\cr 5^3{:}(4\times A_5).2\cr 7^2{:}(3\times 2A_4)\cr5^2{:}2A_5\cr \hline
\end{array}
\;
\begin{array}{c}
Fi_{24}'\cr\hline Fi_{23}\cr 2\udot Fi_{22}{:}2\cr (3\times O_8^+(3){:}3){:}2\cr
O_{10}^-(2)\cr 3^7\udot O_7(3)\cr 3^{1+10}{:}U_5(2){:}2\cr 2^{11}\udot M_{24}\cr
2^2\udot U_6(2){:}S_3\cr 2^{1+12}\udot3\udot U_4(3).2_2\cr
3^{2+4+8}.(A_5\times 2A_4).2\cr 3^3.[3^{10}].GL_3(3)\cr(A_4\times O_8^+(2){:}3){:}2\cr He{:}2 \mbox{ (two)}\cr
2^{3+12}.(L_3(2)\times A_6)\cr 2^{6+8}.(S_3\times A_8)\cr (3^2{:}2\times G_2(3))\udot 2\cr
(A_5\times A_9){:}2\cr A_6\times L_2(8){:}3\cr7{:}6\times A_7\cr U_3(3){:}2 \mbox{ (two)}\cr
 L_2(13){:}2 \mbox{ (two)}\cr 29{:}14\cr \hline
\end{array}
\;
\begin{array}{c}
Fi_{24}\cr\hline Fi_{24}'\cr Fi_{23}\times 2\cr (2 \times 2\udot Fi_{22}){:}2\cr
S_3\times O_8^+(3){:}S_3\cr O_{10}^-(2){:}2\cr 3^7\udot O_7(3){:}2\cr
3^{1+10}{:}(U_5(2){:}2\times 2)\cr 2^{12}\udot M_{24}\cr
(2 \times 2^2\udot U_6(2)){:}S_3\cr 2^{1+12}\udot3\udot U_4(3).2^2\cr
3^{2+4+8}.(S_5\times 2S_4)\cr 3^3.[3^{10}].(L_3(3)\times 2^2)\cr S_4\times O_8^+(2){:}S_3\cr
2^{3+12}.(L_3(2)\times S_6)\cr
2^{7+8}.(S_3\times A_8)\cr (S_3\times S_3\times G_2(3)){:}2\cr
S_5\times S_9\cr 
S_6\times L_2(8){:}3\cr7{:}6\times S_7\cr 7^{1+2}{:}(6\times S_3).2\cr 29{:}28\cr \hline
\end{array}
$$

$$
\begin{array}{c}
\mathbb M \mbox{ (first part)}\cr\hline
2\udot B\cr 2^{1+24}\udot Co_1\cr 3\udot Fi_{24}\cr 2^2\udot {}^2E_6(2){:}S_3\cr
2^{10+16}\udot O_{10}^+(2)\cr 2^{2+11+22}\udot(M_{24}\times S_3)\cr
3^{1+12}\udot 2\udot Suz{:}2\cr 2^{5+10+20}\udot(S_3\times L_5(2))\cr
S_3\times Th\cr 2^{3+6+12+18}\udot (L_3(2)\times 3S_6)\cr 3^8\udot O_8^-(3).2_3\cr
(D_{10}\times HN)\udot 2\cr (3^2{:}2\times O_8^+(3))\udot S_4\cr
3^{2+5+10}.(M_{11}\times 2S_4)\cr 3^{3+2+6+6}{:}(L_3(3)\times SD_{16})\cr
\hline
\end{array}
\;
\begin{array}{c}
\mathbb M \mbox{ (second part)}\cr\hline
5^{1+6}{:}2\udot J_2{:}4\cr (7{:}3\times He){:}2\cr (A_5\times A_{12}){:}2\cr
5^{3+3}\udot(2\times L_3(5))\cr (A_6\times A_6\times A_6).(2\times S_4)\cr
(A_5\times U_3(8){:}3_1){:}2\cr 5^{2+2+4}{:}(S_3\times GL_2(5))\cr
(L_3(2)\times S_4(4){:}2)\udot 2\cr 7^{1+4}{:}(3\times 2S_7)\cr
(5^2{:}[2^4]\times U_3(5)).S_3\cr (L_2(11)\times M_{12}){:}2\cr
(A_7\times (A_5\times A_5){:}2^2){:}2\cr 5^4{:}(3 \times 2\udot L_2(25)){:}2_2\cr
7^{2+1+2}{:}GL_2(7)\cr M_{11}\times A_6\udot 2^2\cr\hline
\end{array}
\;
\begin{array}{c}
\mathbb M \mbox{ (third part)}\cr\hline
(S_5\times S_5\times S_5){:}S_3\cr (L_2(11)\times L_2(11){:}4\cr 13^2{:}2L_2(13).4\cr
(7^2{:}(3\times 2A_4)\times L_2(7)).2\cr (13{:}6\times L_3(3))\udot 2\cr
13^{1+2}{:}(3\times 4S_4)\cr L_2(71)\cr L_2(59)\cr 11^2{:}(5\times 2A_5)\cr
L_2(41)\cr
L_2(29){:}2\cr 7^2{:}SL_2(7)\cr L_2(19){:}2\cr 41{:}40 \cr \mbox{others?}\cr\hline
\end{array}
$$


\end{document}